\documentclass[a4paper, 11pt]{article}
\usepackage[top=0.8 in,bottom=0.8 in,left=0.8 in,right=0.8 in]{geometry}
\usepackage{amsmath,amsthm}
\usepackage{amssymb}
\usepackage[colorlinks,linkcolor=blue,citecolor=blue]{hyperref}
\usepackage{float}
\usepackage{cite}
\usepackage{graphicx,subfigure}
\usepackage{subfigure}
\usepackage{epic}
\usepackage[all]{xy}
\usepackage{tikz}
\usepackage{color}
\usetikzlibrary{intersections}
\usetikzlibrary{cd}
\usetikzlibrary{positioning}
\usepackage{multirow} 
\usepackage{cprotect}
\usepackage{hyperref}

\newtheorem{lemma}{Lemma}[section]
\newtheorem{theorem}{Theorem}[section]

\newtheorem{remark}{Remark}[section]

\def \mc#1{\mathcal{#1}}

\def \mb#1{\mathbb{#1}}
\def \mf#1{\mathbf{#1}}
\def \mr#1{\mathrm{#1}}
\def \op#1{\operatorname{#1}}

\title{Solving the Quispel-Roberts-Thompson maps using Kajiwara-Noumi-Yamada's representation of elliptic curves}

\author{Xing Li$^{1,2}$ and  Tomoyuki Takenawa$^2$\\
${}^1$
  Department of Mathematics, Shanghai University\\
  Shanghai 200444, P.R. China\\
  ${}^2$
  Tokyo University of Marine Science and Technology\\
  2-1-6 Etchu-jima, Koto-Ku, Tokyo, 135-8533, Japan\\
email: xing26@shu.edu.cn \quad
takenawa@kaiyodai.ac.jp\\
}

\date{}

\begin{document}

\maketitle

\begin{abstract}
It is well known that the dynamical system determined by a Quispel-Roberts-Thompson map (a QRT map) preserves a pencil of biquadratic polynomial curves on ${\mathbb{CP}}^1 \times {\mathbb{CP}}^1$.  In most cases this pencil is elliptic, i.e. its generic member is a smooth algebraic curve of genus  one, and  the system can be solved as a translation on the elliptic fiber to which the initial point belongs. However, this procedure is rather complicated to handle, especially in the normalization process. In this paper, for a given initial point on an invariant elliptic curve, we present a method to construct the solution directly in terms of  the Weierstrass sigma function, using  Kajiwara-Noumi-Yamada's parametric representation of elliptic curves.

\vskip 6pt

\end{abstract}

\section{Introduction}

The QRT maps are a family of two-dimensional maps introduced by Quispel, Roberts and Thompson \cite{QRT89} which preserve a one-parameter family of biquadratic curves on ${\mathbb{CP}}^1 \times {\mathbb{CP}}^1$:
\begin{equation}\label{invariantK}
P(x,y; K)=\mathbf{x^T Ay}+K \mathbf{x^T B y}=0
\end{equation}
and a measure \cite{RQ92}
\begin{equation}\label{measure}
\omega=\frac{dx \wedge dy}{P(x,y; K)},
\end{equation}
where  $\mathbf{x}=(x^2,x,1)^T$, $\mathbf{y}=(y^2,y,1)^T$, $\mathbf{A}=(a_{ij})_{0\leq i,j\leq2}$ and $\mathbf{B}=(b_{ij})_{0\leq i,j\leq2}$ are coefficient matrices in $\op{Mat}_{3\times 3}(\mathbb{C})$, and $K$ stands for the single parameter. Such a one-parameter family of curves called a pencil as it depends linearly on the parameter $K$. Moreover, since a horizontal (or vertical) line intersects the generic biquadratic pencil only at two points, it is natural to define a horizontal switch $r_x: (x,y)\to (\bar{x},y)$ and a vertical switch $r_y: (x,y)\to (x, \bar{y})$ on the biquadratic curve. Thus a QRT map can be decomposed into two involutions as
$\varphi=r_y\circ r_x$ \cite{IR02, Tsuda04}:
\begin{equation}\label{qrt}
\bar{x}=\frac{f_1(y)-f_2(y)x}{f_2(y)-f_3(y)x}, \qquad
\bar{y}=\frac{g_1(\bar{x})-g_2(\bar{x})y}{g_2(\bar{x})-g_3(\bar{x})y},
\end{equation}
where $\mathbf{f}=(f_1(y),f_2(y),f_3(y))^T=\mathbf{Ay \times B y}$ and $\mathbf{g}=(g_1(\bar{x}),g_2(\bar{x}),g_3(\bar{x}))^T=\mathbf{A^T \bar{x} \times B^T \bar{x}}$.

Since it is a biquadratic curve in ${\mathbb{CP}}^1 \times {\mathbb{CP}}^1$, each member of the pencil  \eqref{invariantK} is an elliptic curve if it is smooth, otherwise it is a (possibly reducible) rational curve.
In highly degenerate cases, the pencil could be rational, i.e. all the curves are rational, while
in generic cases the pencil is elliptic, i.e. most members are elliptic while finitely many members are rational.
In this paper we focus on the latter case with a given initial point on an  elliptic curve in the pencil.
In this case the dynamical system obtained by iterating QRT map can be solved as a translation on the elliptic curve.
For example, a smooth biquadratic curve can be transformed into the  Weierstrass standard form, but complicated variable transformations are required, and the translation depends on the addition formula of the elliptic curve \cite{Tsuda04}.
Alternatively, if the biquadratic is symmetric, the elliptic curve can be transformed into the canonical form and the solution can be written with Jacobi elliptic functions \cite{book-Baxter82, Veselov-1991}.
In \cite{RCGO-2002} Ramani and his collaborators considered the elliptic parametrization of asymmetric biquadratic curve, and obtained the same standard form as the symmetric case. This result is natural because all elliptic curves can be transformed to the  Weierstrass standard form, but it requires a long variable transformation by computer algebra.
The complex normal form of biquadratic curve has been  given independently in \cite{IR02} with appropriate modular transformations. When the variables are real, the classification and corresponding elliptic parametrization have been given in \cite{IR02, Iatrou-2003}.

Since QRT maps are geometrically simple objects that also appear in useful physical phenomena, they have been studied from various viewpoints, including dynamical properties and properties as subgroups of the automorphism group of rational elliptic surfaces \cite{book-Duistermaat}.

Another important aspect of QRT maps is their relation to the discrete Painlev\'e equations: by de-autonomization of QRT maps, the discrete Painlev\'e equations have been constructed \cite{GRP91, RGH91, Sakai01, CDT2017}.
In the process of this research, various geometric methods have been developed, such as the use of Picard groups \cite{Sakai01} and the construction of embedding maps of elliptic torus in ${\mathbb{CP}}^1 \times {\mathbb{CP}}^1$.

In this paper, we present a method to construct the solution directly in terms of the Weierstrass sigma function, using an embedding of the elliptic torus proposed by Kajiwara, Noumi and Yamada in \cite{KMNOY03}, the pullback of the Picard group of elliptic surface by the embedding \cite{ET2005, CDT2017}, and the elliptic integration. In addition, we give a practical method for elliptic integrals, in which the integration paths must be considered on a double covering. Although this method is classical, especially when the coefficients are complex, controlling the branch cuts so that the square root function is continuous along the path on the covering is complicated and also an obstacle to calculate elliptic integrals numerically, therefore, we provide a Mathematica notebook to achieve this process.

This paper is organized as follows. In Section \ref{sec-2}, as preliminaries, we introduce the initial conditions space of the QRT map, the Picard group, and an embedding of an elliptic torus in ${\mathbb{CP}}^1 \times {\mathbb{CP}}^1$. In Section \ref{sec-3} we present a method to construct the solution using the Weierstrass sigma function. In Section \ref{sec-4}, we apply our method to a specific example.
A Mathematica notebook for computing the example in Section \ref{sec-4} is provided as supplemental material.

\section{Preliminaries}\label{sec-2}
The indeterminate points of QRT map \eqref{qrt}, where the numerator and denominator are simultaneously zero, can be obtained by solving
\begin{equation*}
\left(f_2(y)^2 -f_1(y)f_3(y)=0 \mbox{ and } x=\frac{f_1(y)}{f_2(y)}\right)
\mbox{ or } \left(g_2(\bar{x})^2 -g_1(\bar{x})g_3(\bar{x})=0 \mbox{ and } y=\frac{g_1(\bar{x})}{g_2(\bar{x})}\right)
\end{equation*}
which also give the base points of the pencil.
 Generically,  the elliptic pencil has eight base points counted with multiplicities, which can be removed by blowing up. By resolving these base points, the  obtained space of initial conditions of the QRT map is a rational elliptic surface.

Let $ \mc X$ be a rational surface obtained by blowing up at arbitrary points $p_k$, $k=1,\cdots,K$, possibly infinitely near.
A divisor on $ \mc X$ is a finite formal linear combination $D=\sum_i n_i D_i$, where $D_i$'s are prime divisors, $n_i\in \mb Z$, and the degree of the divisor $D$ is given by the sum of its coefficients, i.e. $\sum_i n_i$. The set of all divisors is an Abelian group denoted as Div$(\mc X)$. Especially, in each blow-up, the point $p_i$ is replaced by one $\mathbb{CP}^1$-line, and hence an exceptional divisor $E_i$ is added.
We say the divisors, $D$ and $D'$ are linearly equivalent if and only if $D-D'$ is a principal divisor that is a divisor of a rational function on $\mc X$. Using this linear equivalence relation,  we can define quotient group $\mr{Cl}(\mathcal{X})=\mr{Div}(\mathcal{X})/\equiv$, where the class of the divisor $D$ is denoted as  $\mc D$.
Since $\mc X$ is a nonsingular irreducible variety, the Picard group Pic($\mathcal{X}$) is isomorphic to $\mr{Cl}(\mathcal{X})$. Let $\mathcal{E}_i$ be the classes of total transform of $E_i$ with respect to blowing up, and by $\mc H_x$ and $\mc H_y$ we denote the classes of total transforms of any vertical axes and horizontal axes. Then the Picard group is a lattice generated by them as
\[
\mr{Pic}(\mc X)=\mr{Span}_{\mb Z}\{\mc H_x,\mc H_y,\mc E_1,\cdots,\mc E_{K}\}.
\]

For the $K$-family of biquadratic curves, when we fix $K$ as $K_0$, $P(x,y; K_0)$ is a biquadratic polynomial on  $\mathbb{CP}^1\times \mathbb{CP}^1$, the closure of all zeros of $P(x,y;K_0)$: $\{(x,y)|P(x,y;K_0)=0\}$ is an algebraic curve, denoted by $\Gamma$. We also use $P(x,y;K_0)=0$ to stand for an algebraic curve.
  In what follows, in order to describe the dynamic behavior of a QRT map in a more direct way, we use an embedding of the elliptic torus $ \mb T$ in $\mb{CP}^1\times \mb{CP}^1$ with elliptic functions of order $2$
\begin{subequations}\label{embedding-1}
\begin{equation*}
\iota: \mb T=\mb C/\mb L \to \mb {CP}^1\times \mb {CP}^1, \qquad u \mapsto \iota(u)=(x,y),
\end{equation*}
\begin{align}
&(x,y)=\left(c\, \frac{[u-\alpha][u-\beta]}{[u-\gamma][u-\delta]},\,
c'\, \frac{[u-\alpha'][u-\beta']}{[u-\gamma'][u-\delta']}\right), \label{KNY_iota}\\
& \alpha+\beta=\gamma+\delta, \quad \alpha'+\beta'=\gamma'+\delta', \label{KNY_abcd}
\end{align}
\end{subequations}
where $\mb L$ is the lattice with $(w_1,w_2)$ as a basic generator: $\mb L=\mb{Z} w_1+\mb{Z} w_2$ and
$[u]$ denotes the Weierstrass sigma function $\sigma(u|\mb L)=\sigma(u;g_2,g_3)$, or the odd theta function $\vartheta(u|\frac{w_2}{w_1})$.

\begin{remark}\label{remark-1} \cite{WW27}
The Weierstrass sigma function and the odd Jacobi theta function are both holomorphic functions and can be used to construct elliptic functions. Below we briefly introduce their properties and the relationship:
\begin{itemize}
\item[(i)]
A lattice $\mb L$ can be generated by more than one pair of generators, such as $(w_1, w_2)$, and $(\widetilde{w}_1,\widetilde{w}_2)$, subject to
\begin{equation*}
\widetilde{w}_1=a w_1+b w_2,\quad \widetilde{w}_2=c w_1+d w_2,\quad  |ad-bc|=1,
\end{equation*}
but the two pairs give rise to unique invariants $g_2$ and $g_3$. Conversely, starting from given invariants $g_2$ and $g_3$ with $g_2^3-27g_3^2\neq0$, the lattice $\mb L$ can be uniquely determined. Therefore we can identify the sigma function  as $\sigma(u)=\sigma(u|\mb L)=\sigma(u;g_2,g_3)$ with $g_2$ and $g_3$ as given invariants.
\item[(ii)]
The Weierstrass sigma function
\begin{equation}
\sigma(u|\mb L) = u \prod_{ (m,n) \in \mb Z^2 \atop (m,n)\neq (0,0) }
\left( 1 - \frac{u}{L_{m,n}} \right) \exp \left( \frac{u}{L_{m,n}} + \frac{u^2}{2L_{m,n}^2} \right),
\end{equation}
 $L_{m,n}= m w_1 +n w_2$, is odd and quasi-periodic
\begin{equation}\label{sigma}
\sigma(u+w_j)=-e^{2 \eta_j(u+w_j/2)}\sigma(u),\quad j=1,2,
\end{equation}
where $\eta_j=\zeta(w_j/2)$ and $\zeta(u)= d (\ln\sigma(u))/du$.
\item[(iii)]
The odd theta function
\begin{equation}
\vartheta\left(u\bigg|\frac{w_2}{w_1}\right)= 2 \sum_{n=0}^{\infty}(-1)^nq^{(n+\frac{1}{2})^2}\sin((2n+1)z),\quad q=\exp\left(\frac{\pi i w_2}{w_1}\right),
\end{equation}
is also quasi-periodic
\begin{align*}
\vartheta\left(u+\pi\bigg|\frac{w_2}{w_1}\right)=-\vartheta\left(u\bigg|\frac{w_2}{w_1}\right),\quad
\vartheta\left(u+\pi \frac{w_2}{w_1}\bigg|\frac{w_2}{w_1}\right)=-q^{-1}e^{-2iu}\vartheta\left(u\bigg|\frac{w_2}{w_1}\right).
\end{align*}
\item[(iv)]
The connection between $\sigma(u|\mb L)$ and $\vartheta(z|\frac{w_2}{w_1})$ is given by
\begin{equation}
\sigma(u|\mb L)=\frac{w_1}{\pi}\exp\left(\frac{\eta_1 u^2}{w_1}\right)\frac{1}{2}q^{-\frac{1}{4}}\prod_{n=1}^{\infty}\left(1-q^{2n}\right)^{-3}\vartheta\left(\frac{\pi u}{w_1}\bigg|\frac{w_2}{w_1}\right),\quad q=\exp\left(\frac{\pi i w_2}{w_1}\right).
\end{equation}
\end{itemize}
\end{remark}
 In this paper, we mainly discuss the case when $[u]$ is the sigma function, and the case of the odd theta function can be obtained similarly. For a fixed lattice $\mb L$, we can obtain an elliptic torus $\mb T$ that is homeomorphic to the quotient group $\mb C/\mb L$.
  From (ii) of the above remark, we immediately obtain that the embedding $\iota$ \eqref{embedding-1} leads to an elliptic parametrization of $(x, y)$, so we have the following Lemma.
\begin{lemma}\label{lemma-1}
The function
\begin{equation}
F(u)=F(u|\mb L; \alpha, \beta, \gamma, \delta) = \frac{[u-\alpha][u-\beta]}{[u-\gamma][u-\delta]}
\end{equation}
with $\alpha+\beta=\gamma+\delta$ satisfies $F(u+w_i)=F(u)$ for $i=1,2$ and $ F( \lambda u | \lambda \mb L;  \lambda \alpha,  \lambda \beta,  \lambda \gamma,  \lambda \delta) = F(u | \mb L; \alpha, \beta, \gamma, \delta) $ for $\lambda \in \mb C$.
\end{lemma}

\section{Solving QRT maps}\label{sec-3}
In this section, we consider the initial value problem for QRT maps with a given initial point on an invariant elliptic curve, and get the solution step by step.

For a given QRT map $\varphi=r_y\circ r_x$ with initial value $(x_0,y_0)$ in \eqref{qrt}, the invariant $K$ can be fixed as $K_0=-\frac{\mathbf{x_0^T B y_0}}{\mathbf{x_0^T Ay_0}}$ and the corresponding curve $\Gamma=\Gamma(K_0)$ is preserved.
Further, we will determine the values of parameters  in \eqref{KNY_iota} from $\Gamma$,
 which are given a priori in constructing the elliptic Painlev\'e equation in \cite{KNY17}. When  $\Gamma$ is singular, a case-by-case study is needed with replacing the elliptic integration by Sakai's period map \cite{Sakai01}.\\

\noindent {\it 0. Check the smoothness of the biquadratic curve}\\
For a biquadratic polynomial $P(x,y)=\sum_{i,j=0}^2 a_{ij} x^i y^j$, the partial discriminant of the $P(x,y)$ with respect to $y$
is a quartic polynomial of $x$:
\begin{equation}
\Delta(x)=\mathop{\rm Discriminant}_y (P(x,y))
= c_0 +c_1 x+c_2 x^2+c_3 x^3+c_4 x^4,
\end{equation}
and $\mathop{\rm Discriminant}(\Delta(x))$ denotes the discriminant of quartic polynomial. Let us introduce the Eisenstein invariant \cite[Section 2.3.5]{book-Duistermaat}
\begin{equation}
\Delta := \frac{\mathop{\rm Discriminant}(\Delta(x))}{256}
= g_2^3 -27 g_3^2,
\end{equation}
where
\begin{align*}
&g_2=c_0 c_4-\frac{c_1 c_3}{4}+\frac{c_2^2}{12},\\
&g_3= -\left(\frac{c_0 c_3^2}{16}+\frac{c_1^2 c_4}{16}-\frac{c_0 c_2 c_4}{6}-\frac{c_1 c_2 c_3}{48}+\frac{c_2^3}{216}\right).
\end{align*}
Then, $\Delta$ is a homogeneous polynomial of degree $12$ in $a_{ij}$, while it coincides with the discriminant
obtained through the Weierstrass standard form as in \cite[Appendix C]{CDT2017}.
It is well known that the corresponding curve (the closure of $P(x,y)=0$ in $\mathbb{CP}^1 \times \mathbb{CP}^1$) is smooth if and only if  $\Delta \neq 0$.

If $\Delta \neq 0$ for the invariant $P(x,y;K)=0$, go to the next step, otherwise stop the procedure.\\

\noindent {\it 1. Move  two points on the elliptic curve to $(\infty, \infty)$ and $(0,0)$}\\
This procedure is not necessary\footnote{
If this step is skipped, $F_{12}(u)$ and $G_{12}(u)$ of \eqref{embedding-2b} is modified as 
\begin{equation*}
F_{12}(u)=\frac{[u-e_{x_2}][u-h_x+e_{x_2}]}{[u-e_{x_1}][u-h_x+e_{x_1}]},\quad G_{12}(u)=\frac{[u-e_{y_2}][u-h_y+e_{y_2}]}{[u-e_{y_1}][u-h_y+e_{y_1}]},
\end{equation*}
where $e_{x_1}$, $e_{x_2}$, $e_{y_1}$ and $e_{y_2}$ correspond to the intersection points of the curve with the lines 
$x=\infty$, $x=0$, $y=\infty$ and $y=0$ respectively.
}, but it would reduce confusion in computing process.
Choose any two distinct points $p_1:(x,y)=(a_1,b_1)$, $p_2:(x,y)=(a_2,b_2)$ on the elliptic curve ($p_1$ and $p_2$ could be the base points, but not necessarily so). 
With the change of variables
\begin{equation}\label{trans-xy}
\rho(x,y)=(\tilde{x},  \tilde{y})=\left(\frac{x-a_2}{x-a_1}, \quad \frac{y-b_2}{y-b_1}\right),
\end{equation}
 we still write $\tilde{x}$, $\tilde{y}$ as $x, y$, the elliptic curve $\Gamma$  passes through the points $ p_1=(\infty,\infty)$, $p_2=(0,0)$, and can be written as follows:
\begin{equation}\label{Gammaxy}
P(x,y;K_0)=a_{21} x^2y+ a_{12}xy^2+a_{20}x^2+a_{11}xy+a_{02}y^2+a_{10}x+a_{01}y=0.
\end{equation}
The initial value $(x_0,y_0)$ is converted to $\rho(x_0, y_0)$.\\

\noindent {\it 2. Periods $w_1$ and $w_2$}\\
By virtue of the symplectic form \eqref{measure},
the periods of elliptic curve $\Gamma$ can be obtained
\cprotect\footnote{On Mathematica, using the invariants $g_2$ and $g_3$, one can also compute the periods  by the command:\\
\verb| 2 WeierstrassHalfPeriods[{g_2,g_3}]|,
 the difference from the periods obtained by the elliptic integral is a scalar multiplication, and the embedding \eqref{KNY_iota} needs to be adjusted by Lemma \ref{lemma-1}.}. Consider $P(x,y;K_0)=0$ as an equation in the variable $y$, the square root of partial discriminant
\begin{align*}
\sqrt{\bigtriangleup(x)}&=\sqrt{-4 x (a_{02} + a_{12} x) (a_{10} + a_{20} x) + (a_{01} + a_{11}x + a_{21} x^2)^2}\\
&=\sqrt{a_{21}^2(x-q_1)(x-q_2)(x-q_3)(x-q_4)}
\end{align*}
has four different branch points at $x=q_i$ $(1\leq i\leq 4)$. By $[q_1,q_2]$ and $[q_3,q_4]$ we denote simple paths form $q_1$ to $q_2$ and from $q_3$ to $q_4$ in $\mb{CP}^1$ without common points, then $[q_1,q_2]$ and $[q_3,q_4]$ can be defined as branch cuts of the Riemann sphere with two  Riemann sheets $\Omega_{\pm}$. For a general $x$, there are two values of $y$, determined by
\begin{equation*}
y_{\pm}= \frac{- (a_{01} + a_{11}x + a_{21} x^2)\pm \sqrt{\Delta(x)}}{2 (a_{02} + a_{12} x)}.
\end{equation*}
The residue of $\omega$ at $y_{+}$ gives rise to $1$-form on the elliptic curve
\begin{align*}
\underset{y=y_{+}}{\mathrm{Res}} \,\omega = \frac{1}{2 \pi i} \oint_{y=y_{+}} \omega=\frac{dx}{\sqrt{\Delta(x)}} ,
\end{align*}
then the periods are given by the value of the integrals
\begin{align*}
w_1=\int_{\delta_1} \underset{y=y_{+}}{\mathrm{Res}}\,\omega,\quad w_2=\int_{\delta_2} \underset{y=y_{+}}{\mathrm{Res}}\,\omega,
\end{align*}
where the path $\delta_1$ is a suitable simple closed curve enclosing the points $q_1$, $q_2$ and the path $\delta_2$ can be taken as a suitable simple closed curve enclosing the points $q_2$, $q_3$.
Thus, we obtain an elliptic torus $\mb T =\mb C/ (\mb Z w_1+\mb Zw_2) $.
Note that if we normalize 1-form on the elliptic curve, i.e. periods $w_1$ and $w_2$ are divided by a constant, the embedding \eqref{KNY_iota} is changed only multiplitically as the second assertion in Lemma~\ref{lemma-1}. \\

\noindent
{\it 3. Selecting $\alpha, \beta, \gamma, \delta, \alpha', \beta', \gamma', \delta'$}\\
In order to obtain the parameters in embedding \eqref{embedding-1}, we introduce the Picard group of the elliptic torus $\mb T$, and the pull-back map $\iota^*$ on the level of the Picard groups. The Picard group of the elliptic torus $\mb T$ consists of the set of all finite formal
linear combinations $D=\sum_i n_i t_i$, where $t_i \in \mb T$ and $n_i\in \mb Z$, modulo principal divisors of elliptic functions ($\equiv$). Let $T_n$ be the set of all classes of divisors of degree $n$, then the Picard group can be decomposed to mutually disjoint sets as
\begin{equation*}
\operatorname{Pic}(\mb T)=\bigcup_{n=-\infty}^{n=+\infty}T_n,
\end{equation*}
where each $T_n$ is isomorphic to $\mb T$.
Let us consider $\iota$ as an embedding from
$\mb T=\mb C/\mb L$ to $\mc X$, where $\mc X$ is a rational surface obtained by blowing up from ${\mathbb{CP}}^1 \times {\mathbb{CP}}^1$ at $p_1$ and $p_2$\footnote{The following argument is still valid if  one changes the definition of $\mc X$ so that $\mc X$ is obtained by blowing up from the corresponding rational elliptic surface at $p_1$ and $p_2$.}. Then
the pull-back map $\iota^*: \op{ Pic}(\mc X)\to \op{Pic}(\mb T)$ gives rise to
\begin{equation}
h_x = \iota^*(\mc H_x),\quad h_y = \iota^*(\mc H_y), \quad e_i = \iota^*(\mc E_i),\quad i=1,2,
\end{equation}
where $h_x, h_y \in T_2$ and $e_i \in T_1$. Especially $\iota^*(\mc E_2) = \iota^*((0,0))$ indicates that one of $\alpha$, $\beta$ is $e_2$ and one of $\alpha'$, $\beta'$ is also $e_2$. Similarly considering $\iota^*(\mc E_1) = \iota^*((\infty, \infty))$, we can assume that
\begin{equation}
\alpha =\alpha' =e_2, \quad \gamma =\gamma' =e_1.
\end{equation}
Moreover, since $h_x$ is the pull-back of the divisor class $\mc H_x$, we have
\begin{equation}
\alpha +\beta =\gamma + \delta = h_x,
\end{equation}
and hence
\begin{equation}
\beta =h_x- e_2, \quad \delta = h_x -e_1.
\end{equation}
In the same way, we have
\begin{equation}
\beta' =h_y- e_2, \quad \delta' = h_y -e_1.
\end{equation}

In summary, as shown in \cite{CDT2017}, the embedding \eqref{embedding-1} can be written  as
\begin{subequations}\label{embedding-2}
\begin{equation} \label{embedding-2a}
\iota(u)= (c_1 \,F_{12}(u),\, c_2\, G_{12} (u)),
\end{equation}
with
\begin{equation}\label{embedding-2b}
F_{12}(u)=\frac{[u-e_2][u-h_x+e_2]}{[u-e_1][u-h_x+e_1]},\quad G_{12}(u)=\frac{[u-e_2][u-h_y+e_2]}{[u-e_1][u-h_y+e_1]},
\end{equation}
\end{subequations}
which was originally proposed by Kajiwara, Noumi and Yamada in \cite{KNY17}.\\

\noindent
{\it 4. Practical method to compute values of $h_x$, $h_y$, $e_i$ and $c_i$}\\
We should start by fixing an arbitrary point $p^O=(x^O, y^O)$ on $\Gamma$, then the inverse point $\iota^{-1}(p^O)=O$ can be labeled as the origin on the elliptic torus $\mb T$.
In the following we formulate a practical method to compute the values of $h_x$, $h_y$ and $e_i$ by using elliptic integral, similar to the step calculating periods $w_1$ and $w_2$.

Recall that $e_1$, $h_x-e_1$, $h_y-e_1$, $e_2$, $h_x-e_2$ and $h_y-e_2$ correspond to the points
\begin{equation*}
\begin{array}{lll}
e_1: p_1=(\infty,\infty), \quad  &h_x-e_1: (\infty, y_1'), \quad &h_y-e_1: (x_1', \infty),\\
e_2: p_2=(0,0), \quad &h_x-e_2: (0, y_2'), \quad &h_y-e_2: (x_2', 0),\\
\end{array}
\end{equation*}
where  $x_2'$ is determined by solving $P(x_2',0; K_0)=0$ with $x_2'\neq 0$ and
$x_1'$ is determined in chart-$(x,y=1/Y)$
\begin{equation*}
Y^2 P(x_1', 1/Y; K_0) |_{Y=0}=0 \quad \mbox{with $x_1'\neq 0$}
\end{equation*}
and so on.
Hence for example $e_1$ can be computed by the integration
\begin{align}
e_1=\int_{x^O}^{\infty} \underset{y=y_{+}}{\mathrm{Res}}\,\omega,
\end{align}
where the terminal point of the integral path corresponds to the point $(\infty, \infty)$.
But it should be careful that we have two Riemann sheets $\Omega_{\pm}$ with infinity, so in order to distinguish two sheets, a valid way is to consider the integration paths on double covering, such as $e_1$ and $h_x-e_1$ are provided by
\begin{align}
e_1=\int_{(x^O,y^O)}^{(\infty, \infty)} \underset{y=y_{+}}{\mathrm{Res}}\,\omega, \qquad
h_x-e_1=\int_{(x^O,y^O)}^{(\infty, y_1')} \underset{y=y_{+}}{\mathrm{Res}}\,\omega,
\end{align}
the second component of the target point, i.e. the value of $y$, indicates that the Riemann sheet where $(\infty, y_1')$ is located is different from the Riemann sheet where $(\infty, \infty)$  is located.
Hence, in practical computation, we need to control the branch cuts of square root so that the value of $y$ changes continuously along the path on the covering, and $h_y-e_1$, $e_2$, $h_x-e_2$, $h_y-e_2$ should be taken into account in the same way.

Finally, coefficients $c_1$ and $c_2$ in \eqref{embedding-2a} can be determined by the pull-back $\iota^*((x_1', \infty))$ (or $\iota^*((x_2',0))$) and  $\iota^*(( \infty, y_1'))$ (or $\iota^*((0,y_2'))$) respectively
\begin{equation*}
\begin{split}
c_1&=\frac{x_1'}{F_{12}(h_y-e_1)}=\frac{x_2'}{F_{12}(h_y-e_2)},\\
c_2&=\frac{y_1'}{G_{12}(h_x-e_1)}=\frac{y_2'}{G_{12}(h_x-e_2)}.
\end{split}
\end{equation*}
In light of above calculation, the parameters in embedding $\iota$ \eqref{embedding-2} are uniquely determined.\\

\begin{remark} \label{ambiguity}
When one chooses different paths for integration, the obtained value could become $e_1'=e_1 + n_{1} w_1 + m_{1} w_2$,
similarly,  the values of  $e_2$ and $h_x-e_1$ could become $e_2'=e_2 + n_{2} w_1 + m_{2} w_2$,
and $h_x'-e_1'=h_x-e_1 + n_{3} w_1 + m_{3} w_2$ ($n_{i},{m_{i}} \in \mb Z$ for all $i$, $j$).
According to (ii) of Remark~\ref{remark-1}, the more general formula is
\begin{equation}
\sigma(u +  n_1 w_1 + m_1w_2) = c \cdot \exp\left( 2(n_1 \eta_1+m_1 \eta_2) u\right) \sigma(u),
\end{equation}
where $c$ is a non-zero constant, therefore there exists a non-zero constant $c'$ such that
\begin{align}
\frac{[u-e_2'][u-h_x'+e_2']}{[u-e_1'][u-h_x'+e_1']} = c'  \frac{[u-e_2][u-h_x+e_2]}{[u-e_1][u-h_x+e_1]}
\end{align}
holds ( in exponential part, the terms containing $u$ are eliminated, since the sum of zero points of the numerator and that of the denominator are the same).
\end{remark}

\begin{figure}[h!]
\centering
\includegraphics[width=0.4\linewidth]{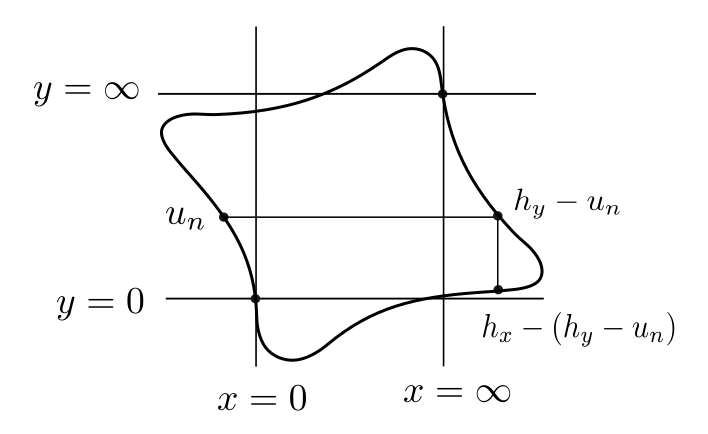}
\caption{Dynamics of the solution of QRT map $\widetilde{\varphi}=\rho \circ r_y\circ r_x \circ \rho^{-1}$}
\label{fig1}
\end{figure}

\noindent {\it 5. Action of QRT maps on the elliptic torus}\\
Recall the geometric description of QRT map $\varphi=r_y\circ r_x$ \eqref{qrt},
 by $\mathcal X$ we still denote the rational surface under coordinate transformation $\rho$ \eqref{trans-xy},
then the QRT map $\widetilde{\varphi}=\rho \circ r_y\circ  r_x \circ\rho^{-1}: \mathcal X\to \mathcal X$.
Let us consider the dynamic behavior on the elliptic torus  $\mathbb T$, we can get  $\iota^{-1}\circ \rho \circ r_x\circ\rho^{-1}\circ \iota$: $u \mapsto h_y-u$, $\iota^{-1}\circ \rho \circ r_y\circ\rho^{-1}\circ \iota$ acts as $u \mapsto h_x-u$, and the dynamics can be simply viewed as the Figure \ref{fig1}.
We immediately get the transformation $\tau(u)=\bar{u}=u-(h_y-h_x)$, which gives the commutative diagram
\begin{figure}[H]
\begin{tikzpicture}
\begin{tikzcd}[column sep=60pt,row sep=35pt]
& \mc X \arrow[r, "\widetilde{\varphi}"] & \mc X \\
&\mb T\arrow[u, "\iota"]\arrow[r, "{\tau}"] & \mb T \arrow[u, swap, "\iota"]
\end{tikzcd}
\end{tikzpicture}
\hspace{0.5em}
\begin{tikzpicture}
\begin{tikzcd}[column sep=50pt,row sep=30pt]
& \op{Pic}(\mc X) \arrow[d, "\iota^*"] & \op{Pic}(\mc X) \arrow[d, swap, "\iota^*"]\arrow[l, swap, "\widetilde{\varphi}^*"] \\
&\op{Pic}(\mb T) & \op{Pic}(\mb T)\arrow[l, "{\tau^*}"]
\end{tikzcd}
\end{tikzpicture}
\end{figure}
Combining it with embedding \eqref{embedding-2}, we obtain the main result.
\begin{theorem}\label{th}
Rewrite QRT map \eqref{qrt} as
\begin{equation}
x_{n+1}=\frac{f_1(y_n)-f_2(y_n)x_n}{f_2(y_n)-f_3(y_n)x_n}, \qquad
y_{n+1}=\frac{g_1(x_{n+1})-g_2(x_{n+1})y_n}{g_2(x_{n+1})-g_3(x_{n+1})y_n},
\end{equation}
with initial value $(x_0,y_0)$, then the solution of the QRT map corresponding to the initial value is
 \begin{equation}\label{sol}
(x_n,y_n)=\rho^{-1}\big(c_1 F_{12}(u_0+n(h_x-h_y)),~c_2 G_{12}(u_0+n(h_x-h_y))\big),
\end{equation}
where  $\rho(x,y)$ is defined as in \eqref{trans-xy},  $F_{12}(u)$, $G_{12}(u)$ are defined as in \eqref{embedding-2b} and $u_0$ is determined by the integral from $(x^O, y^O)$ to  $\rho(x_0, y_0)$ as
\begin{align}
u_0=\int_{(x^O,y^O)}^{\rho(x_0, y_0)} \underset{y=y_{+}}{\mathrm{Res}}\,\omega.
\end{align}
\end{theorem}

\section{Example}\label{sec-4}
In this section, we will show more details through an example.

Let $\varphi_1$ be the QRT map defined by the matrices
\begin{subequations}
\begin{equation}
\mathbf{A}=\left(\begin{array}{ccc}
0 &-7-i & 3+i \\
4 i& -5+2 i&2-i\\
3+4 i&6 & 0
\end{array}\right),\quad
 \mathbf{B}=\left(\begin{array}{ccc}
0&0&0\\
0&0&1\\
0&1&0
\end{array}\right).
\end{equation}
\end{subequations}
Then $\varphi_1=r_y\circ r_x$ is written as
\begin{align}
x_{n+1}&=\frac{ y_n\left(1+i\right) \left(-(7+i) x_n y_n+(3+i) x_n+4 i y_n^2-(8+2 i) y_n-4-i\right)}{2\left((3+4 i) y_n-1-2 i\right) \left(x_n+y_n\right)},
\\
y_{n+1}&=-\frac{i x_{n+1} \left(-4 i x_{n+1} y_n+(7+i) x_{n+1}^2+(8-i) x_{n+1}-(3+4 i) y_n-4-i\right)}{\left(4 x_{n+1}+4-3 i\right) \left(x_{n+1}+y_n\right)}.
\end{align}

Take the initial value $(x_0, y_0)=(1, 0.437561\, +0.328195 i)$, then it is easy to find $\mf{x^T A y}=0$, and so $K$ is $0$. If the subscript $n$ is omitted, $ P(x_n,y_n;0)=P(x,y;0)$ is written as
\begin{equation}\label{Gamma}
\begin{split}
&P(x,y;0)=
-(7+i) x^2 y+4 i x y^2+(3+i) x^2-(5-2 i) x y+(3+4 i) y^2+(2-i) x+6 y,
\end{split}
\end{equation}
and the corresponding elliptic curve $\Gamma$ is the closure of $\{(x,y)|P(x,y;0)=0\}$. Considering $P(x,y;0)=0$ as an equation with respect to $y$, we get four branch points of the square root of the discriminant $\sqrt{\Delta(x)}$:
\begin{align*}
  &\{q_1,\, q_2,\, q_3,\, q_4\}\\
=&\{-1.69314+0.647424 i,-1.01244+0.358514 i,0.264181\, -0.0620125 i,1.083\, +0.827275 i\}.
\end{align*}
Take the $1$-form on the elliptic curve $\Gamma$ as ${\mathrm{Res}}_{y=y_{+}}\,\omega=dx/\sqrt{\Delta(x)}$ and select two simple closed curves $\delta_1$ and $\delta_2$ enclosing branch points $q_1, q_2$ and $q_2, q_3$ respectively. Then we get the periods
\begin{align*}
&w_1=\int_{\delta_1} \underset{y=y_{+}}{\mathrm{Res}}\,\omega=0.0207773\, +0.438853 i,\\
&w_2=\int_{\delta_2}  \underset{y=y_{+}}{\mathrm{Res}}\,\omega=
-0.59573+0.114127 i,
\end{align*}
and the elliptic torus $\mb T =\mb C/ (\mb Z w_1+\mb Zw_2)$.

\begin{figure}[!h]
\centering
\begin{minipage}{5cm}
\includegraphics[width=\textwidth]{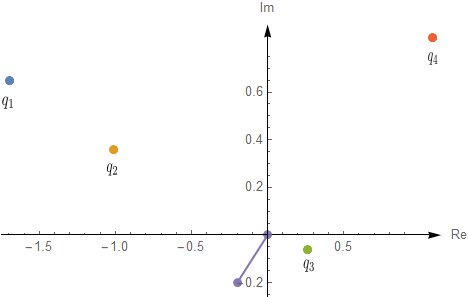}
\centering
\footnotesize{(a) $e_2$}
\end{minipage}
\hspace{0in}
\begin{minipage}{5cm}
\includegraphics[width=\textwidth]{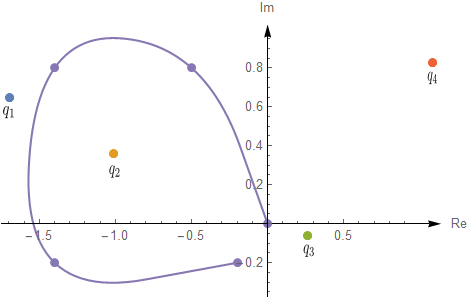}
\centering
\footnotesize{(b)  $h_x-e_2$}
\end{minipage}
\begin{minipage}{5cm}
\includegraphics[width=\textwidth]{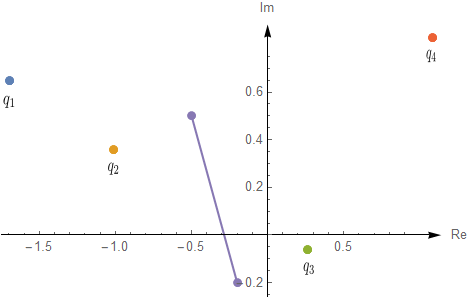}
\centering
\footnotesize{(c)  $h_y-e_2$}
\end{minipage}
\hspace{0in}
\begin{minipage}{5cm}
\includegraphics[width=\textwidth]{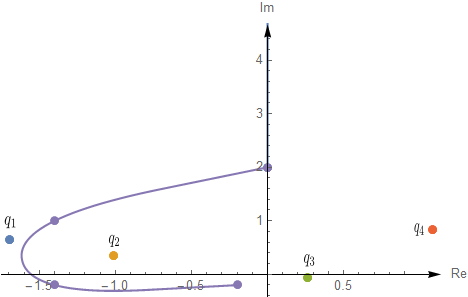}
\centering
\footnotesize{(d) $e_1$}
\end{minipage}
\begin{minipage}{5cm}
\includegraphics[width=\textwidth]{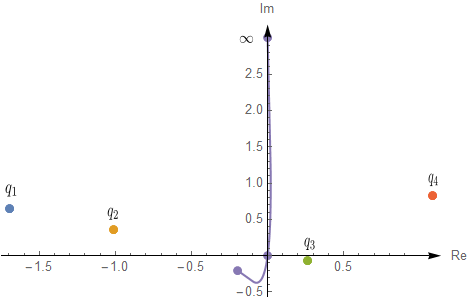}
\centering
\footnotesize{(e) $h_x-e_1$}
\end{minipage}
\begin{minipage}{5cm}
\includegraphics[width=\textwidth]{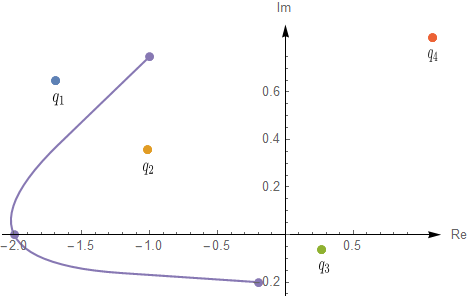}
\centering
\footnotesize{(f) $h_y-e_1$}
\end{minipage}
\caption{The integration paths}
\label{paths}
\end{figure}

Taking $(x^O,y^O)=(-0.2-0.2 i,0.0864885\, -0.00825559 i)$ as the starting point of the integration path on $\Gamma$ and solving $P(x,0;0)=0$ and $P(0,y;0)=0$, we obtain $x_2'=-0.5+0.5 i$ and $y_2'=-0.72+0.96 i$. Then the integration paths shown in Figure \ref{paths} give rise to
\begin{align*}
  e_2&=\int_{(x^O,y^O)}^{(0, 0)} \underset{y=y_{+}}{\mathrm{Res}}\,\omega
 =0.0302102\, +0.0268586 i,\\
 h_x- e_2&=\int_{(x^O,y^O)}^{(0, y_2')} \underset{y=y_{+}}{\mathrm{Res}}\,\omega
 =-0.406377+0.0917916 i,\\
  h_y- e_2&=\int_{(x^O,y^O)}^{(x_2', 0)} \underset{y=y_{+}}{\mathrm{Res}}\,\omega
=-0.0486106+0.092694 i.
\end{align*}
The value of $e_1$ can be computed along the same way.
 Let us take terminal points as $(x_1',\infty)=(-1.+0.75 i,\infty)$ and $(\infty, y_1')=(\infty, 0.44\, +0.08 i)$  by solving $Y^2 P(x, 1/Y; 0) |_{Y=0}=0$ and $X^2 P(1/X, y; 0) |_{X=0}$ $=0$ respectively and take the corresponding integration paths in Figure \ref{paths}, we have:
\begin{align*}
  e_1&=\int_{(x^O,y^O)}^{(\infty, \infty)} \underset{y=y_{+}}{\mathrm{Res}}\,\omega
 =-0.367314-0.171699 i,\\
 h_x'- e_1&=\int_{(x^O,y^O)}^{(\infty, y_1')} \underset{y=y_{+}}{\mathrm{Res}}\,\omega
 =-0.0085321+0.290344 i,\\
 h_y'- e_1&=\int_{(x^O,y^O)}^{(x_1', \infty)} \underset{y=y_{+}}{\mathrm{Res}}\,\omega
=-0.267552-0.0334266 i.
\end{align*}
Here we noticed that $e_2+(h_x-e_2)=e_1+(h_x'-e_1)$, but $e_2+(h_y-e_2)=e_1+(h_y'-e_1)+w_1-w_2$.
However, as discussed in Remark~\ref{ambiguity}, the only difference is the constant multiplication, and  we can take $h_y$ as $h_y=e_2+(h_y-e_2)$.

Finally, the coefficients can be determined as
\begin{equation*}
\begin{split}
c_1&=\frac{x_1'}{F_{12}(h_y-e_1)}=\frac{x_2'}{F_{12}(h_y-e_2)}
=0.673966\, -2.34158 i,\\
c_2&=\frac{y_1'}{G_{12}(h_x-e_1)}=\frac{y_2'}{G_{12}(h_x-e_2)}
=3.67383\, +4.1975 i.
\end{split}
\end{equation*}
As a result, we obtain the solution
\begin{equation} \label{sol-2}
\begin{aligned}
x_n&=c_1 F_{12}(u_0+n(h_x-h_y)), \qquad y_{n}&=c_2 G_{12}(u_0+n(h_x-h_y)),
\end{aligned}
\end{equation}
where $F_{12}(u)$, $G_{12}(u)$ are defined as in \eqref{embedding-2b}, and $u_0$ is given by
\begin{align*}
  u_0&=\int_{(x^O,y^O)}^{(x_0,y_0)} \underset{y=y_{+}}{\mathrm{Res}}\,\omega
 =0.0714408\, -0.115967 i.
\end{align*}
\begin{remark}
The solution of a QRT map with a given initial point is determined only by the invariant biquadratic to which the initial point belongs (not by the whole pencil of the invariant curves).
For example, if we construct a new QRT map $\varphi_2$ by keeping matrix $\mathbf{A}$ and choosing a different matrix $\mathbf{B}'$ as
\begin{equation}
\mathbf{A}=\left(\begin{array}{ccc}
0 &-7-i & 3+i \\
4 i& -5+2 i&2-i\\
3+4 i&6 & 0
\end{array}\right),\quad
 \mathbf{B'}=\left(\begin{array}{ccc}
0&0&0\\
0&1&0\\
0&0&0
\end{array}\right),
\end{equation}
then the QRT map $\varphi_2$ is
\begin{align}
x_{n+1}&=-\frac{y_n \left(1+ i\right) \left(6+(3+4 i) y_n\right)}{2 x_n \left((3+4 i) y_n-1-2 i\right)},\\
y_{n+1}&=-\frac{ x_{n+1}(1+2 i) \left(1+(1+i) x_{n+1}\right)}{\left(4 x_{n+1}+4-3 i\right) y_n},
\end{align}
and the invariant K-family on ${\mathbb{CP}}^1 \times {\mathbb{CP}}^1$ is:
\begin{equation}
 P'(x,y;K)= \mathbf{x^T Ay}+K \mathbf{x^T B' y}=0.
\end{equation}
If we further use the same initial value $(x_0,y_0)=(1, 0.437561\, +0.328195 i)$ so that $K=0$, and $P'(x,y;0)$ is the same as \eqref{Gamma}, then the QRT map $\varphi_2$ has the same solution \eqref{sol-2} with $\varphi_1$.
\end{remark}

\section*{Acknowledgements}
In December 2019, TT was asked by Prof. Saburo Kakei whether the QRT dynamics could be solved directly using knowledge of the initial value space for the Painlev\'e equations. Although TT could not give a good answer at the time, the inspiration for this paper came from that conversation, and the authors deeply grateful to him. The authors thank Prof. Da-jun Zhang and Prof. Anton Zhamay for giving XL an opportunity to study at TUMSAT with TT.  The authors also thank the anonymous reviewers for their valuable comments on the manuscript.
TT was supported by the Japan Society for the Promotion of Science, Grand-in-Aid
(C) (17K05271). XL was supported by the NSF of China (11875040) and Science and Technology Innovation Plan of Shanghai (20590742900).

\end{document}